\documentclass{agtart_a}
\pdfoutput=1

\usepackage{graphicx}

\title{Quantum link invariant from\\the Lie superalgebra~${\mathfrak D}_{2\,1,\alpha}$}
\author{Bertrand Patureau-Mirand}
\givenname{Bertrand}
\surname{Patureau-Mirand}
\address{LMAM Universit{\'e} de Bretagne-Sud\\
Centre de Recherche\\Campus de Tohannic\\\newline
BP 573\\F-56017 Vannes\\France}
\email{bertrand.patureau@univ-ubs.fr}
\urladdr{http://www.univ-ubs.fr/lmam/patureau/}



\received{1 February 2005}
\revised{}
\accepted{15 August 2005}
\proposed{}
\seconded{}
\publishedonline{12 March 2006}
\published{12 March 2006}

\volumenumber{6}
\issuenumber{}
\publicationyear{2006}
\papernumber{11}
\startpage{329}
\endpage{349}

\doi{}
\MR{}
\Zbl{}

\arxivreference{math.GT/0404548}

\AtBeginDocument{\let\bar\wbar\let\tilde\wtilde\let\hat\what}
\newcommand{\we}{\smash{\rlap{\kern 6pt 
\raise 4pt\hbox{\footnotesize $\sim$}}}\longrightarrow}

\newcommand{\wb}{\overline}
\newcommand{\wbtot}{\widetilde}
\newcommand{\wt}{\widetilde}

\newcommand{\Mod}{\operatorname{Mod}}
\newcommand{\Obj}{\operatorname{Obj}}

\renewcommand{\Z}{Z_{\operatorname{ad}}}
\newcommand{ \ad}{{\operatorname{ad}}}
\newcommand{\Id}{\operatorname{Id}}
\newcommand{\e}{\operatorname{e}}
\newcommand{\LL}{{\cal{L}}}
\newcommand{\T}{{\cal{T}}}
\newcommand{\Tp}{{\cal{T}}_p}
\newcommand{\A}{{\cal{A}}}

\newcommand{\calD}{{\cal{D}}}
\newcommand{\ZZ}{\mathbb{Z}}  
\newcommand{\QQ}{\mathbb{Q}}
\newcommand{\RR}{\mathbb{R}}
\newcommand{\CC}{\mathbb{C}}
\newcommand{\NN}{\mathbb{N}}
\renewcommand{\L}{\Lambda}
\renewcommand{\S}{\mathfrak{S}}

\newcommand{\go}{\longrightarrow}
\newcommand{\im}{\mapsto}

\newcommand{\osp}{\mathfrak{osp}}
\newcommand{\DD}{{{\mathfrak D}_{2\,1}}}
\newcommand{\DDa}{{{\mathfrak D}_{2\,1,\alpha}}}
\newcommand{\Uh}{{{U_h}\DD}}
\newcommand{\sll}{\mathfrak{sl}}

\renewcommand{\sp}{\mathfrak{sp}}
\newcommand{\Xd}{\chi_{\DD}}
\newcommand{\ext}{\wedge}

\renewcommand{\Q}[2]{{[#1]_#2}}
\newcommand{\et}{\quad\textrm{and}\quad}
\newcommand{\ets}{\quad,\quad}

\makeatletter
\def\cnewtheorem#1[#2]#3{\newtheorem{#1}{#3}[section]
\expandafter\let\csname c@#1\endcsname\c@prop}

\theoremstyle{plain}
\newtheorem{prop}{Proposition}[section]
\cnewtheorem{lem}[prop]{Lemma}
\cnewtheorem{theo}[prop]{Theorem}
\cnewtheorem{coro}[prop]{Corollary}
\cnewtheorem{conj}[prop]{Conjecture}
\theoremstyle{definition}
\cnewtheorem{Rq}[prop]{Remark}
\cnewtheorem{defi}[prop]{Definition}

\makeatother  %  move after \newtheorem block

\newcommand{\epsw}[2]
         {\begin{array}{c} \hspace{-.3mm}
        \raisebox{-4pt}{\includegraphics[width=#2]{\figdir/#1}}
        \hspace{-.6mm}\end{array}}
\newcommand{\epsh}[2]
         {\begin{array}{c} \hspace{-.3mm}
        \raisebox{-4pt}{\includegraphics[height=#2]{\figdir/#1}}
        \hspace{-.6mm}\end{array}}
\newcommand{\pstext}[1]{\hspace{1mm}\raisebox{-0.5ex}{%
\includegraphics[height=2.5ex]{\figdir/#1}}\hspace{1mm}}
\newcommand{\psdiag}[3]{\hspace{1mm}\raisebox{-#1mm}{%
\includegraphics[height=#2mm]{\figdir/#3}}\hspace{1mm}}

\begin{document}

\begin{asciiabstract}
  The usual construction of link invariants from quantum groups
  applied to the superalgebra D_{2 1,alpha} is shown to be
  trivial. One can modify this construction to get a two variable
  invariant. Unusually, this invariant is additive with respect to
  connected sum or disjoint union. This invariant contains an infinity
  of Vassiliev invariants that are not seen by the quantum invariants
  coming from Lie algebras (so neither by the colored HOMFLY-PT nor by
  the colored Kauffman polynomials).
\end{asciiabstract}

\begin{htmlabstract}
  The usual construction of link invariants from quantum groups applied
  to the superalgebra D<sub>2&nbsp;1,&alpha;</sub> is shown to be
  trivial. One can modify this construction to get a two variable
  invariant. Unusually, this invariant is additive with respect to
  connected sum or disjoint union. This invariant contains an infinity
  of Vassiliev invariants that are not seen by the quantum invariants
  coming from Lie algebras (so neither by the colored HOMFLY-PT nor by
  the colored Kauffman polynomials).
\end{htmlabstract}

\begin{abstract}
  The usual construction of link invariants from quantum groups applied
  to the superalgebra ${\mathfrak D}_{2\,1,\alpha}$ is shown to be
  trivial. One can modify this construction to get a two variable
  invariant. Unusually, this invariant is additive with respect to
  connected sum or disjoint union. This invariant contains an infinity
  of Vassiliev invariants that are not seen by the quantum invariants
  coming from Lie algebras (so neither by the colored HOMFLY-PT nor by
  the colored Kauffman polynomials).
\end{abstract}

\maketitle

%%%%%%%%%%%%%%%%%%%%%%%%%%%%%%%%%%%%%%%%%%%%%%%%%%%%%%%%%%%%%%%%%%%%%%%%%%%%%%
%%%%%%%%%%%%%%%%%%%%%%%%%%%%%%%%%%%%%%%%%%%%%%%%%%%%%%%%%%%%%%%%%%%%%%%%%%%%%%
\section*{Introduction}
In his classification of finite dimensional Lie superalgebras \cite{Kac},
V\,G Kac introduces a family of simple Lie
superalgebras $\DDa$ depending of the parameter
$\alpha\in\CC\setminus\{0,-1\}$. The notation evokes a deformation of
the Lie superalgebra $\osp(4,2)$ which is obtained for
$\alpha\in\{-2,-\frac12,1\}$.

There is a method to construct framed link invariants with a
deformation of an enveloping Lie algebras. It follows from work of
Drinfel'd \cite{Dr} that for a fixed simple Lie algebra, all
deformations give the same link invariant. This is not
clear{\footnote{At the time of writing this paper; but see
Geer \cite{Geer1,Geer2} for new results.}} for the simple Lie
superalgebra $\DDa$. So we explore here two possibly different
deformations of $U\DDa$: the Kontsevich--Drinfel'd deformation and
the quantum group $\Uh$ described by Y\,M Zou and H Thys
\cite{Zou,Th}. The two corresponding link invariants will
be denoted by the letters $Z$ and $Q$.

First we will see that the quantum link invariants $Z_{\DDa,V}$
obtained from any representation $V$ of the Lie superalgebra $\DDa$ is
determined by the linking matrix (see \fullref{ZD} and 
\fullref{linking}). Similarly the quantum link invariants $Q_{\DDa,L}$
obtained from the adjoint representation of $\DDa$ is constant equal
to $1$. (See \fullref{QD}).

A similar problem was encountered by J Murakami \cite{Mur93}, Kashaev
\cite{Kv} and Degushi \cite{D}: the quantum invariants they considered
factor by the zero quantum dimension (the invariant of the
unknot). The remedy is to ``divide'' by this quantum dimension by
considering $(1,1)$--tangles instead of links. Here the invariant of
the planar trivalent tangle $\Theta$ is zero and we give a
construction to ``divide'' trivalent tangles and links by this
$\Theta$.

From this we construct a map $\wt Z$ in \fullref{GZD} that associates to
each framed link an element of the ring
$\QQ[[a_1,a_2,a_3]]^{\S_3}_{/(a_1+a_2+a_3)}$ (quotient of the ring of
symmetric series in three variables) and in \fullref{GQD} we construct
$\wbtot Q$ that associates to
each framed link an element of the ring
$$\ZZ[\frac12,q_1,q_2,q_3, \Q41^{-1},\Q42^{-1},\Q43^{-1}]
/(q_1q_2q_3=1).$$

The existence and the invariance of $\wt Z$ is easy to proof but one
can hardly compute it. On the other side, it takes much more work to
proof that $\wbtot Q$ is well defined but it can be computed with an
$R$--matrix.

It is natural to conjecture that the two deformations of $U\DDa$ are
equivalent. Knowing this would implies that the two maps $\wbtot Q$
and $\wt Z$ would essentially be the same (setting
$q_i=\e^{\frac{a_i}2}$ etc) and then their value would be in the
intersection of these two different rings.

The author thanks Y\,M Zou for sending his papers and C Blanchet for
reading the first version of this manuscript.  The author also wishes
to thank the referee for numerous helpful comments.

\section{Statement of the results}
We work with framed trivalent tangles and knotted trivalent graphs which are
generalizations of framed tangles and links (they are embeddings of
$1$-$3$--valent graphs in $S^3$). Here are an example of a trivalent tangle
and of a knotted trivalent graph (see \fullref{3net} for precise definition).
$$\epsh{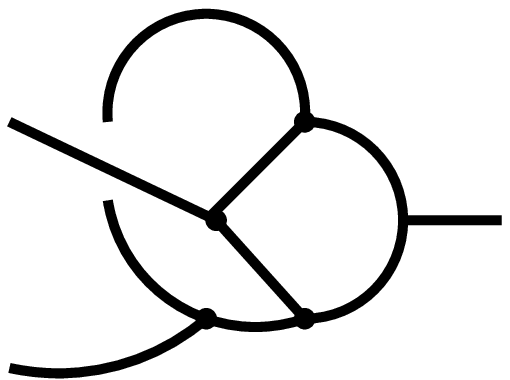}{10ex}\qquad\qquad\epsh{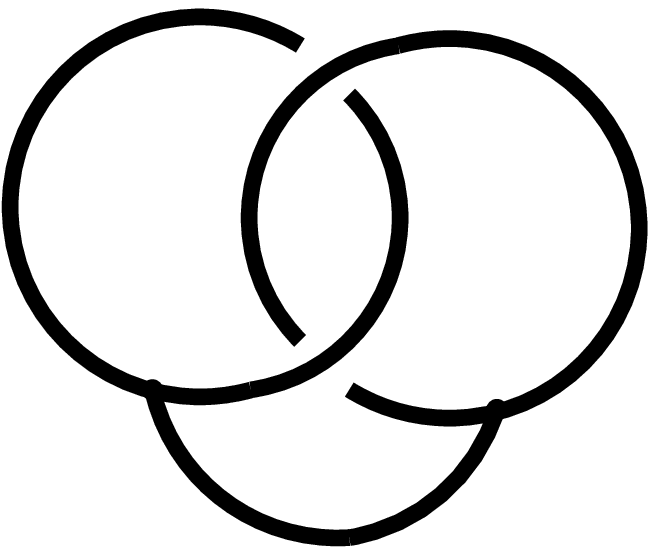}{10ex}$$ 
% %
We will call a framed knotted trivalent graph ``proper'' if it has at least one
trivalent vertex.

The ``adjoint'' Kontsevich integral (cf \fullref{Kad}) associate
to each trivalent tangle a series of $1$-$3$--valent Jacobi
diagrams. When composed with the weight function $\Phi_{\DD}$
associated to the Lie superalgebra $\DD$ (which is a generalization of
$\DDa$), it gives an functor $Z_{\DD,L}$ from the category of
trivalent tangles to a completion of the category of representation of
$\DD$ (here $L$ denotes the adjoint representation of $\DD$).
\begin{prop}\label{Z=0}
 On trivalent tangles, $Z_{\DD,L}$ does not depend of the
  framing. Furthermore, if $N$ is a knotted trivalent graph then
$$Z_{\DD,L}(N)=\begin{cases}
0 & \text{if $N$ is a proper knotted trivalent graph}, \\
1 & \text{if $N$ is a link}
\end{cases}$$
\end{prop}
The ``adjoint'' Kontsevich integral of a knotted trivalent graph lives
in the space of closed $3$--valent Jacobi diagrams. This space has a
summand isomorphic to the algebra $\L$ defined by P Vogel in
\cite{Vo} on which $\Phi_{\DD,L}$ is not trivial. Using this map one
can construct an invariant of knotted trivalent graph $\wt Z$ with
values in the ring $\QQ[[a_1,a_2,a_3]]^{\S_3}_{/(a_1+a_2+a_3)}$.
\begin{prop} $Z_{\DD,L}$ and $\wt Z$ are related by
$$Z_{\DD,L}\left(\epsh{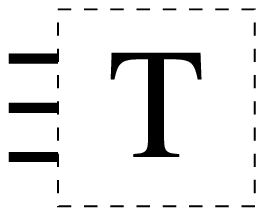}{4ex}\right)= \wt
Z\left(\epsh{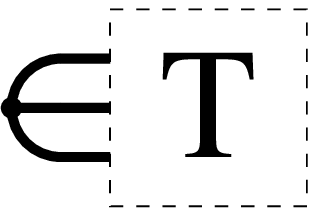}{4ex}\right)
.Z_{\DD,L}\left(\epsh{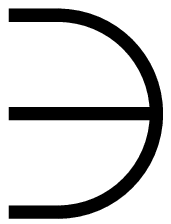}{3ex}\right)$$
for any trivalent tangle $T$ with $3$ $1$--valent vertices.
\end{prop}

We state similar results for $Q_{\DDa,L}$: The quantum group $\Uh$
has an unique topologically free representation $L$ whose classical
limit is the adjoint representation of $\DDa$. This module is autodual
and there is an unique map (up to a scalar) $L\otimes L\go L$ whose
classical limit is the Lie bracket. As usual, coloring a trivalent tangle
with $L$ gives a functor $Q_{\DDa,L}$ from the category of trivalent tangles
to the category of representation of the quantum group $\Uh$ and, in
particular, a knotted trivalent graph invariant.
\begin{prop} \label{Q=0} On trivalent tangles, $Q_{\DDa,L}$ does not depend of
  the framing. Furthermore, if $N$ is a knotted trivalent graph then
$$Q_{\DDa,L}(N)=\begin{cases}
0 & \text{if $N$ is a proper knotted trivalent graph}, \\
1 & \text{if $N$ is a link}
% \text{otherwise}
\end{cases}$$
\end{prop}
We modify this invariant to the following:
\begin{theo} There is an unique invariant of proper knotted trivalent graph $\wt Q$,
  with values in the ring $\ZZ[\frac12,q_1,q_2,q_3,
\Q41^{-1},\Q42^{-1},\Q43^{-1}] /(q_1q_2q_3=1)$ defined by the
property:
$$Q_{\DDa,L}\left(\epsh{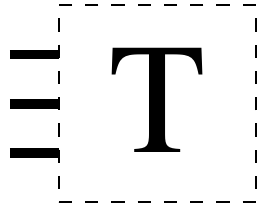}{4ex}\right)= \wt
Q\left(\epsh{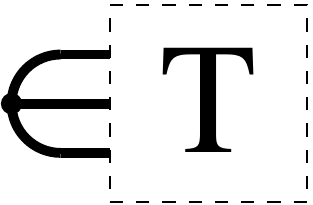}{4ex}\right)
.Q_{\DDa,L}\left(\epsh{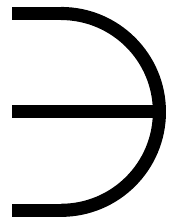}{3ex}\right)$$
for any trivalent tangle $T$ with $3$ univalent vertices.
\end{theo}

\begin{theo}\label{QbKnot}
There exists an invariant of framed links uniquely determined by:
$$\wbtot Q\left(\epsh{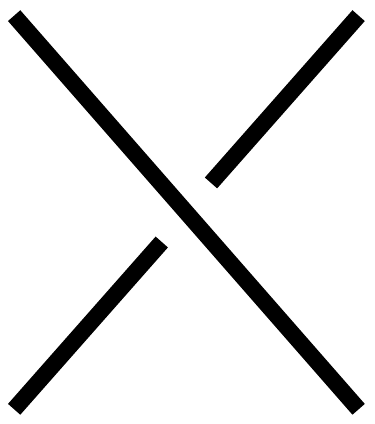}{4ex}\right)-\wbtot
Q\left(\epsh{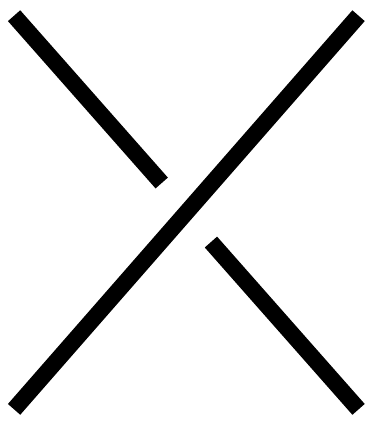}{4ex}\right) =\wt
Q\left(\epsh{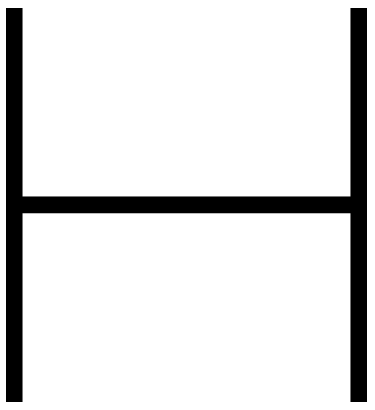}{4ex}\,-\epsh{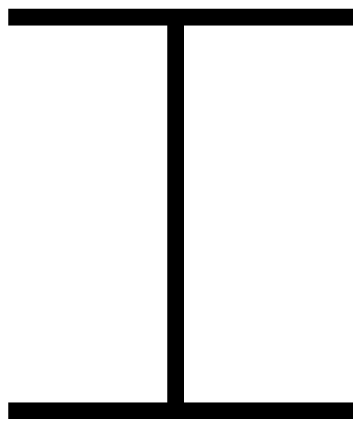}{4ex}+\frac12\left(\epsh{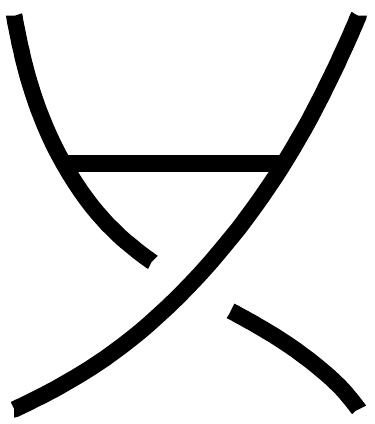}{4ex}+
\epsh{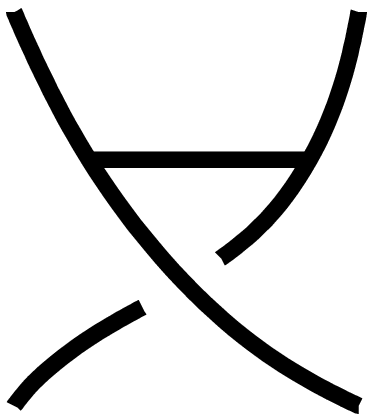}{4ex}\right)\right)$$
$$\text{ and }\qquad \wbtot Q\left(\text{Unlink}\right)=0$$
\end{theo}

\begin{conj} \label{conjI}
$\wbtot Q$  takes value in the polynomial algebra
$\ZZ[\sigma_+,\sigma_-]$
where $\sigma_+=(q_1^2+q_2^2+q_3^2)$ and
  $\sigma_-=(q_1^{-2}+q_2^{-2}+q_3^{-2})$.
\end{conj}
%%%%%%%%%%%%%%%%%%%%%%%%%%%%%%%%%%%%%%%%%%%%%%%%%%%%%%%%%%%%%%%%%%%%%%%%%%%%%%
%%%%%%%%%%%%%%%%%%%%%%%%%%%%%%%%%%%%%%%%%%%%%%%%%%%%%%%%%%%%%%%%%%%%%%%%%%%%%%
\section{The superalgebra $\DDa$}
%%%%%%%%%%%%%%%%%%%%%%%%%%%%%%%%%%%%%%%%%%%%%%%%%%%%%%%%%%%%%%%%%%%%%%%%%%%%%%
\subsection{The Cartan matrix}\label{DDa}
Let $\alpha\in\CC\setminus\{0,-1\}$. The simple Lie superalgebra
$\DDa$ introduced by V\,G Kac \cite{Kac} has the following
Cartan matrix
$$
A_\alpha=(a_{ij})_{\scriptscriptstyle 1\leq i,j\leq 3}=\left(
\begin{array}{ccc}
0 & 1 & \alpha \\
-1 & 2 & 0 \\
-1 & 0 & 2
\end{array}\right),
$$
Where the first simple root is odd and the two others are even. So the
superalgebra is generated by the nine elements $e_i$, $f_i$, $h_i$
($i=1\cdots3$) with the following relations:
$$[h_i,h_j]=0\ets [e_i,f_j]=\delta_{i,j}h_i\ets
[h_i,e_j]=a_{ij}e_j\ets [h_i,f_j]=- a_{ij}f_j\ets$$
$$[e_2,e_3]=[f_2,f_3]=[e_1,e_1]=[f_1,f_1]=[e_i,[e_i,e_1]]=[f_i,[f_i,f_1]]=0\mbox{
  for }i=1,2$$ 
% %
Its even part is isomorphic to the Lie algebra
$L=\sll_2\oplus\sll_2\oplus\sll_2$ and the bracket makes the odd part
an $L$--module isomorphic to the tensor product of the three standard
representations of $\sll_2$. So we can identify the set of weights
with a subset of $\ZZ^3$ such that the three simple roots are:
$$(1,-1,-1)\quad(0,2,0)\quad(0,0,2)$$
Then the set of positive roots is 
$$\{\beta_1=(0,0,2),\,\beta_2=(1,-1,1),\,\beta_3=(1,1,1),\,\beta_4=(2,0,0),$$
$$\qquad\beta_5=(1,-1,-1),\,\beta_6=(1,1,-1),\,\beta_7=(0,2,0)\}$$

%%%%%%%%%%%%%%%%%%%%%%%%%%%%%%%%%%%%%%%%%%%%%%%%%%%%%%%%%%%%%%%%%%%%%%%%%%%%%%
\subsection{The superalgebra $\DD$}\label{DD}
We use in \fullref{weight} the following construction of $\DD$ (see
\cite{Vo} and \cite{Zo2}):

Let $R=\QQ[a_1^{\pm},a_2^{\pm},a_3^{\pm}]_{/(a_1+a_2+a_3)}$ the
quotient of a Laurent polynomial algebra in three variables. Let
$E_1$, $E_2$ and $E_3$ be three two dimensional free $R$--modules
equipped with an non degenerate antisymmetric form:
$\cdot\ext\cdot:\Lambda^2E_i\stackrel{\sim}{\go}R$. We can see $E_i$
as a supermodule concentrated in odd degree equipped with a
supersymmetric form.

The superalgebra $\DD$ is defined as a supermodule by
$$\DD=\sp(E_1)\oplus\sp(E_2)\oplus\sp(E_3)\oplus (E_1\otimes
E_2\otimes E_3)$$
The bracket is defined by the Lie algebra structure on the even part,
by the standard representation of the even part on the odd part and
for the tensor product of two odd elements by the formula:
$$\begin{array}{l}[e_1\otimes e_2\otimes e_3,e'_1\otimes e'_2\otimes e'_3]=\,
\frac12\left(a_1\, e_2\ext e'_2\, e_3\ext
e'_3\,(e_1.e'_1)\hspace{15ex}\ \right.\\
\left.\hspace{15ex}\ \hfill+a_2 \,e_1\ext e'_1 \,e_3\ext e'_3\,(e_2.e'_2) + 
a_3 \,e_1\ext e'_1 \,e_2\ext e'_2\,(e_3.e'_3)\right)\end{array}$$
where $(e_i.e'_i)\in\sp(E_i)$ sends $x\in E_i$ on $e_i\ext x\,e_i'+
e_i'\ext x\,e_i$.

The non degenerate supersymmetric bilinear form on $\DD$ is up to
multiplication by a scalar the orthogonal sum of $-\frac4{a_i}$ times
the killing form of $\sp(E_i)$ plus the tensor product of the three
antisymmetric forms on $E_1\otimes E_2\otimes E_3$.

If $\phi: R\go \CC$ is a ring map, then the complex Lie
superalgebra $\DD\otimes_\phi\CC$ is isomorphic to
$\DDa$ where $\alpha=\frac{\phi(a_3)}{\phi(a_2)}$.

%%%%%%%%%%%%%%%%%%%%%%%%%%%%%%%%%%%%%%%%%%%%%%%%%%%%%%%%%%%%%%%%%%%%%%%%%%%%%%
%%%%%%%%%%%%%%%%%%%%%%%%%%%%%%%%%%%%%%%%%%%%%%%%%%%%%%%%%%%%%%%%%%%%%%%%%%%%%%
\section{The Kontsevich--Drinfel'd invariant}
In the following, if $n\in\NN$ we will denote by $[n]$ the set
$\{1,\,2,\cdots,\,n\}$.

%%%%%%%%%%%%%%%%%%%%%%%%%%%%%%%%%%%%%%%%%%%%%%%%%%%%%%%%%%%%%%%%%%%%%%%%%%%%%%
\subsection{The category of trivalent tangles}\label{3net}
Let $X$ be a finite set. A $X$--diagram is a finite graph, whose
vertices are either $1$--valent or $3$--valent and oriented, with the
data of an isomorphism between $X$ and the set of $1$--valent
vertices. An orientation at a $3$--valent vertex $v$ is a cyclic
ordering on the set of the three edges going to $x$.

Following \cite{BS} we define a trivalent tangle on
$\phi:X\hookrightarrow\RR^3$ as an embedding of an $X$--diagram in
$\RR^3$, with image $N\subset [0,1]\times \RR^2$ together with a
vector field along $N$ such that the points of $N$ lying in the planes
$\{0\}\times \RR^2$ and $\{1\}\times \RR^2$ are exactly
$\phi(X)$. Additionally, we require that the normal plane of $N$ at an
univalent vertex $v$ is parallel to the plane $\{0\}\times\RR^2$, the
vector field assigned to $N$ at $v$ is $(0,0,1)$, and at each
$3$--valent vertex, the orientation of the plane tangent to $N$ given
by the vector field agree with the orientation of the $3$--valent
vertex of the underlying $X$--diagram. When represented by planar
graphs the framing is obtained by taking the vector field pointing
upward.

Two trivalent tangles $T_1$ and $T_2$ are equivalent if one can be
deformed into the other (within the class of trivalent tangles) by a
smooth one parameter family of diffeomorphisms of $\RR^3$.

A framed knotted trivalent graph is a trivalent tangle with no
univalent vertices. We will call a framed knotted trivalent graph
``proper'' if it has at least one trivalent vertex.

Let $M$ be the non-associative monoid freely generated by one letter
noted ``$\circ$''. If $m\in M$, the length of $m$ is the number of
letter in $m$. This gives a partition $M=\bigsqcup_{n\in\NN}M_n$.

The category $\Tp$ (resp.\ $\T$) is the $\QQ$--linear monoidal category
whose set of object is $M$ (resp.\ $\NN$). If $(m,m')\in M_n\times
M_{n'}$ (resp.\ if $(n,n')\in\NN^2$), the set of morphisms
$\Tp(m,m')\simeq \T(n,n')$ is the vector space with bases the set of
equivalence classes of trivalent tangles on the map
($\phi:[n]\sqcup[n']\simeq \left(\{0\}\times[n]\times\{0\}\right)\cup
\left(\{1\}\times[n']\times\{0\}\right)\subset\RR^3$.

The composition is just given by gluing the corresponding univalent
vertices of the tangles. The tensor product of morphism is given by
the juxtaposition of tangles.

%%%%%%%%%%%%%%%%%%%%%%%%%%%%%%%%%%%%%%%%%%%%%%%%%%%%%%%%%%%%%%%%%%%%%%%%%%%%%%
\subsection{The Kontsevich integral for trivalent tangles and the functor $Z_{\DD,L}$}
%%%%%%%%%%%%%%%%%%%%%%%%%%%%%%%%%%%%%%%%%%%%%%%%%%%%%%%%%%%%%%%%%%%%%%%%%%%%%%
\subsubsection{Category of Jacobi diagrams}
We represent an $X$--diagram (or ``Jacobi diagram'') by a graph
immersed in the plane in such a way that the cyclic order at each
vertex is given by the orientation of the plane.

We define the degree of an $X$--diagram to be half the number of the
vertices.

Let $\A(X)$ (resp.\ $\A_c(X)$) denotes the completion (with respect to
the degree) of the quotient of the $\QQ$--vector space with basis the
$X$--diagrams (resp.\ connected $X$--diagrams) by the relations $(AS)$
and $(IHX)$:
\begin{enumerate}
\item If two Jacobi diagrams are the same except for the cyclic order
of one of their vertices, then one is minus the other (relation called
(AS) for antisymmetry).
  $$\lower15pt\hbox{\includegraphics[width=.85cm]{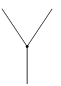}}+
\lower15pt\hbox{\includegraphics[width=.85cm]{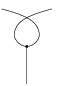}}\equiv0$$
\item The relation (IHX) (or Jacobi) deals with three diagrams which
  differ only in the neighborhood of an edge.
  $$\lower10pt\hbox{\includegraphics[width=.85cm]{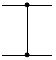}}\equiv\,
\lower10pt\hbox{\includegraphics[angle=90, width=.85cm]{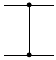}}\,  
-\lower10pt\hbox{\includegraphics[width=.85cm]{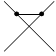}}$$
\end{enumerate}
As we want to work with $\DD$ (which has superdimension 1), we add the
relation that identify the Jacobi diagram with only one circle with
$1$ (this will mean that the superdimension of $\DD$ is $1$). So we
can remove or add some circle to a Jacobi diagram without changing its
value in $\A$.

Let $\calD$ denote the $\QQ$--linear monoidal category defined by
\begin{enumerate}
\item Obj$(\calD)=\left\{[n],n\in\NN\right\}$
\item $\calD([p],[q])=\A([p]\amalg[q])$
\item The composition of a Jacobi diagram from $[p]$ to $[q]$ with a Jacobi diagram
from $[q]$ to $[r]$ is given by gluing the two diagrams along $[q]$.
\item The tensor product of two object is $[p]\otimes [q]=[p+q]$ and
the tensor product of two Jacobi diagrams is given by their disjoint union.
\end{enumerate}

\begin{Rq}
The composition map $\calD([p],[q])\otimes\calD([q],[r])\go\calD([p],[r])$ has
degree $-q$.
\end{Rq}

The algebra $\L$ is the sub-vector space of $\A_c([3])$ formed by
totally antisymmetric elements for the action of the permutation group
$\S_3$.

$\L$ has a natural commutative algebra structure and acts on each
space $\A_c(X)$: If $u$ lies in $\L$ and $K\in \A_c(X)$ is a
$X$--diagram, a Jacobi diagram representative for $u.K$ is obtained by
inserting $u$ at a $3$--valent vertex of $K$. Exceptionally in $\L$,
the degree will be defined by half the number of vertices minus two so
that the unit of $\L$ has degree $0$.

In the following, we will denote by $\A^+$ and $\A_c^+$ (resp.\ by
$\A^0$ and $\A_c^0$) the subspaces generated by Jacobi diagrams having
at least one (resp.\ having no) $3$--valent vertices.

For small $n$ one can describe $\A_c([n])$ (cf \cite{Vo}):

$\A_c([1])$ is zero; $\A_c^+([0])$ and $\A_c^+([2])$ are free
$\L$--modules with rank one, generated by the following elements:
$$\epsh{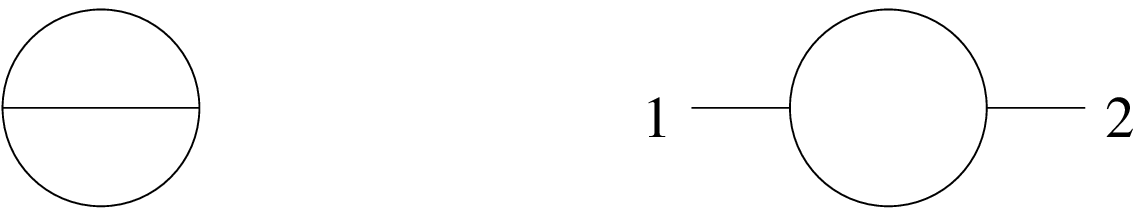}{5ex}$$
Let $\Theta$ be this generator of $\A_c^+(\emptyset)$. Furthermore, we
don't know any counterexample to the following conjecture:
$\A_c([3])=\L$.

$\L$ is generated in degree one by the element $t$:
$$\begin{array}{ccccc}t=&\epsh{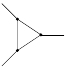}{5ex}&=\frac12&
\epsh{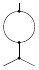}{5ex} \end{array}$$

\subsubsection[The Kontsevich-adjoint functor $Z_{ad}$]{The Kontsevich-adjoint functor $\Z$}\label{Kad}\qua
We follow here A\,B Berger and\break I Stassen \cite[Definition and Theorem
2.8]{BS} who have defined a {\em unoriented universal
Vassiliev--Kontsevich in\-va\-riant} generalized for trivalent tangles
(cf also Murakami and Ohtsuki \cite{MO}). We just consider it for the
``adjoint'' representation so we compose their functor (whose values
are bicolored graph) with the functor that forget the coloring of the
edges.
\begin{theo}[cf \cite{BS}]\label{VK}
There is an unique monoidal functor $\Z : \T_p \go \calD$ (the universal
adjoint Vassiliev--Kontsevich invariant) defined by the following
assignments:

$\Z(m):=[n]$, where $m\in M_n$ is a non-associative word of length
$n$.
\begin{align*}
\Z\Big(\epsh{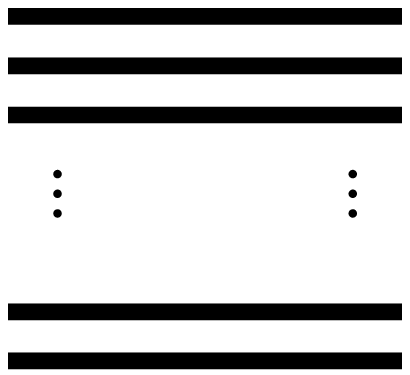}{4ex}:((uv)w)\rightarrow(u(vw))\Big)
\quad&:=\quad \epsh{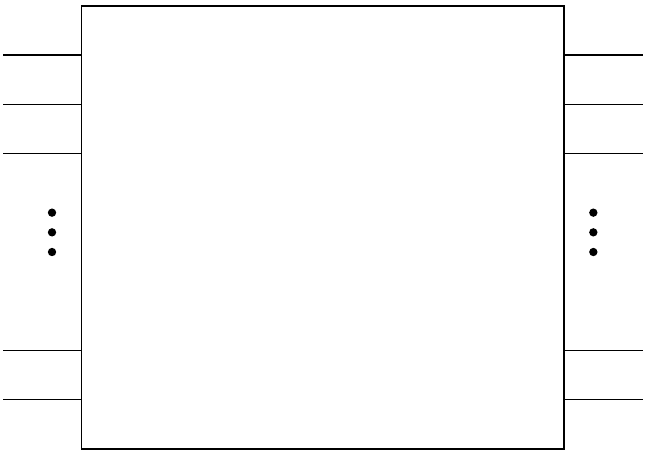}{5ex}\hspace{-29pt}\Phi_{uvw}\\
\Z\left(\epsh{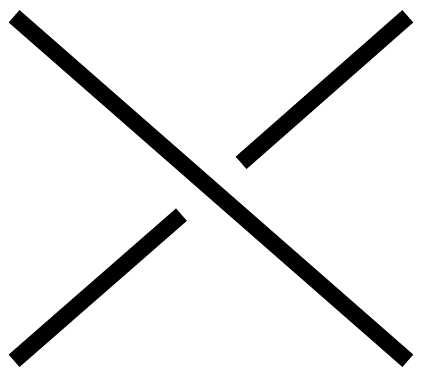}{3ex}\right)\quad&:=\quad\epsh{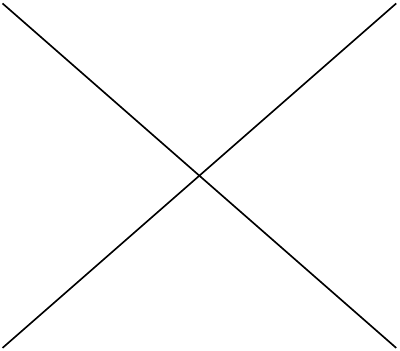}{3ex}\circ
\e^{-\frac{1}{2}\psdiag{0}{2}{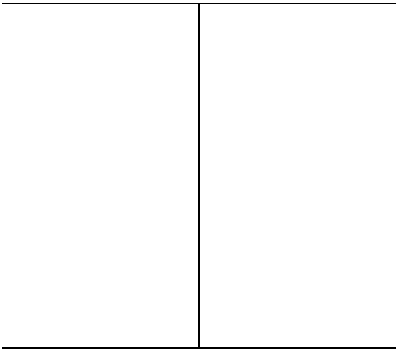}}\\
\Z\left(\epsh{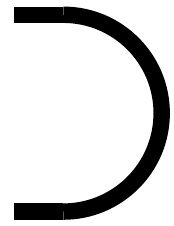}{3ex}\right)\quad&:=\quad\epsh{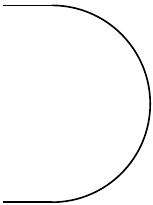}{3ex}\circ
(\Id\otimes C^{\frac12})\\
\Z\left(\epsh{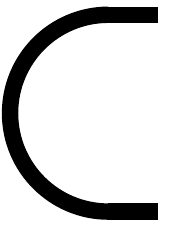}{3ex}\right)\quad&:=\quad(\Id\otimes
C^{\frac12})\circ\epsh{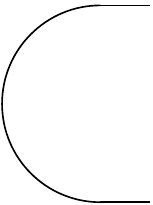}{3ex}\\
\Z\left(\epsh{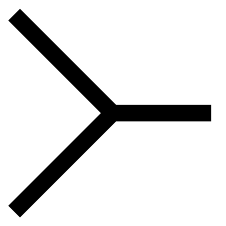}{3ex}\right)\quad&:=\quad
r\epsh{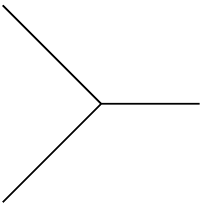}{3ex}\\
\Z\left(\epsh{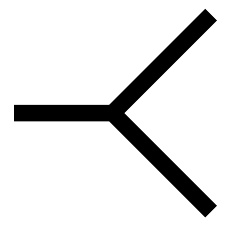}{3ex}\right)\quad&:=\quad
r\epsh{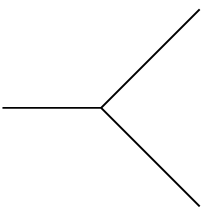}{3ex}\\
\end{align*}
where
\begin{itemize}
\item $e^{\pm \frac{1}{2}\psdiag{0}{2}{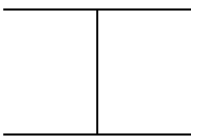}} := \sum_{n=0}^{\infty}$
$(\pm\frac{1}{2})^n\frac{1}{n!}\psdiag{1}{4}{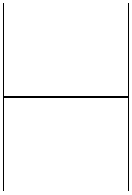}^{\circ n}$ 
\item The elements $\Phi_{uvw}$ are constructed from an even rational
      horizontal associator $\Phi$ with $\Phi^{321}=\Phi^{-1}$ as in
      \cite{LM1}.
\item $C$ :=\pstext{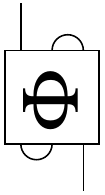}$=Id+\gamma\epsh{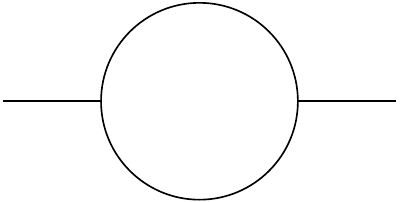}{2ex}$ with
$\gamma\in\L$
\item $r$ can be any element of $\L$. We make the following
  normalization: $r=1$ so that
  $\Z(\Theta\in\Tp(\emptyset,\emptyset))=(1+2t\gamma)\Theta\in\calD(0,0)$
\end{itemize}
\end{theo}
The difference with \cite{BS} for the image of the elementary
trivalent tangles with one $3$--valent vertex is because they use the
Knizhnik--Zamolodchikov associator which has not the good property for
cabling (see \cite{LM1}).

The Kontsevich integral of the unknot has an explicit expression (see
\cite{BNGRT}) but it seems difficult to give an explicit expression of
$\gamma$. Nevertheless it allows to say that $\gamma$ lives in odd
degree and starts with $\gamma=\frac1{24}t+\cdots$.

\subsubsection{Weight functions}\label{weight}
Let $<.,.>_\DD$ denotes the supersymmetric invariant non degenerate
bilinear form on $\DD$) and let $\Omega\in \DD\otimes \DD$ be the
associated Casimir element.

\begin{theo}[cf \cite{Vo}]\label{PhiD}
  There exists an unique $\QQ$--linear monoidal functor $\Phi_{\DD}$
  from $\calD$ to the category $\Mod_\DD$ of representations of $\DD$
  such that:
\begin{enumerate}
\item $\Phi_{\DD}([1])=\DD$ (the adjoint representation).
\item Its values on the elementary morphisms
  $$\put(-105,0){\includegraphics{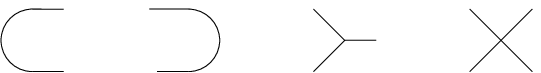}}$$
  are given by:
\begin{enumerate}
\item The Casimir of $\DD$: $\Omega\in \DD^{\otimes 2}\hookrightarrow
  \Mod_\DD(R,(\DD)^{\otimes 2})$
\item The bilinear form $<.,.>_\DD:\DD^{\otimes 2}\go R$
\item The Lie bracket seen as a map in $\Mod_\DD((\DD)^{\otimes
    2},(\DD))$
\item The symmetry operators: $\begin{array}[t]{ccl} \DD^{\otimes
    2}&\go&\DD^{\otimes 2}\\ x\otimes y&\im&(-1)^{|x||y|}y\otimes x \end{array}$
\end{enumerate}
\end{enumerate}
Furthermore, there exists a graded character with value in $\wt
R=\QQ[[\sigma_2,\sigma_3]]$
(here we set $\sigma_2=a_1a_2+a_2a_3+a_3a_1$ and $\sigma_3=a_1a_2a_3$):
$\Xd:\L\go\wt R$ such that: 
$$\forall u\in\L,\, \forall K\in
\A_c^+([p]\amalg[q])\subset\calD([p],[q]),\,
\Phi_{\DD}(uK)=\Xd(u)\Phi_{\DD,\Omega}(K)$$ 
One has $\Xd(t)=0$ and the functor $\Phi_{\DD}$ is zero on the
generators of $\A_c^+(\emptyset)$ and $\A_c^+([2])$.
\end{theo}

\subsubsection{$Z_{\DD,L}$ and the quantum Jacobi relation}\label{ZD}
Composing the adjoint-Kontsevich invariant with the weight function
associated with $\DD$, one get a functor
$Z_{\DD,L}:\T_p\go\Mod_\DD$. For a simple Lie algebra, Drinfel'd
equivalence results for quasi-triangular quasi-Hopf algebra would
imply for the two constructions to give equivalent representations of
$\T_p$ but this is not clear for $\DDa$. So we do the same work for
$Z$:

The functor $Z_{\DD,L}$ produces an invariant of framed knotted trivalent graphs with
values in $\wt R$.

\begin{proof}[Proof of \fullref{Z=0}]
This is a consequence of the fact that $\Phi_\DD$ is zero on any
closed Jacobi diagram having at least one trivalent vertex because $\Phi_\DD$
is $0$ on the generator of $\A_c^+(\emptyset)$.
\end{proof}

\begin{Rq}\label{linking}
This argument can be adapted for other choices of representations of
$\DD$ to prove that the corresponding invariant is determined by the
linking matrix. In particular this is the case for the $\alpha=s$
specialization of the Kauffman polynomial.

$Z_{\DD,L}$ also verify the relations satisfied by $Q_{\DDa,L}$ in
\fullref{relQ}. In particular, as for any simple quadratic Lie
superalgebra, it satisfies the quantum Jacobi relation for
$\phi=-\frac1{2r^2}$ as we will show in a following paper.
\end{Rq}

\subsection{Renormalization of $Z_{\DD,L}$} \label{GZD}
The adjoint Kontsevich integral of a knot $K$ is of the form
$\Z(K)=1+\lambda.\Theta$ for some $\lambda\in\Lambda$. If we apply the
$\DD$ weight system, we just get $1$ since $\Theta$ goes to zero. But
we can ``divide by zero'' defining $\wt Z(K)$ as the weight system
applied to $\lambda$. In the following, we generalize this
construction for links and knotted trivalent graphs.

Let us define $\A'$ to be the quotient of $\A\otimes\L$ by the
relation:

If a Jacobi diagram $K=K_1\sqcup K_2\sqcup\cdots\sqcup
K_n$ represents an element of $\A$ where $K_1$ is a connected Jacobi
diagram such that $K_1=u.K'_1$ for $u\in\L$ then $K\otimes
v=(K'_1\sqcup K_2\sqcup\cdots\sqcup K_n)\otimes u.v$. (This extends
the action of $\L$ to disconnected Jacobi diagrams.)

As before, we define $\calD'$ as the $\L$--linear monoidal category with
the modules $\A'$ as morphisms. 

\begin{prop}
The quotient algebra $\A'(\emptyset)$ is isomorphic to the subalgebra
$\QQ\oplus \Theta\L[\Theta]\subset\L[\Theta]$.

The functor $\Phi_{\DD}$ factor through $p:\calD\go\calD'$.
\end{prop}
\begin{proof}
Just see that $\A(\emptyset)$ is the symmetric algebra on the vector
  space $\A_c(\emptyset)\simeq\QQ\oplus \L\Theta$.

$\Phi_{\DD}$ naturally satisfy the additional relations of $\calD'$ as by
\fullref{PhiD}, it sends via $\Xd$ the elements of $\L$ on scalars.
\end{proof}
We will use the following map on $\A'(\emptyset)$ to get a new
invariant:
\begin{align*}
\Phi'_\DD:\A'(\emptyset)\simeq\QQ\oplus
\Theta\L[\Theta]&\go\QQ[[\sigma_2,\sigma_3]]\\
z+\Theta\lambda+\Theta^2x&\mapsto \Xd(\lambda)
\end{align*}
Where $z\in\QQ$, $\lambda\in\L$ and $x\in\L[\Theta]$.

\begin{defi}For a knotted trivalent graph $L$ set
$$\wt Z(L)=\Phi'_\DD\left(p(\Z(L))\right)$$  
\end{defi}
The planar knotted trivalent graph $\Theta$ is sent by $\Z$ on
$(1+2t\gamma)\Theta\in\A'(\emptyset)$ so
$\wt Z(\Theta)=1$. 

% Let us consider the symmetric bilinear form
% $b:\A'([n])\otimes\A'([n])\go\A'(\emptyset)$ given by gluing the
% diagrams along $[n]$.
Remark that the decomposition $\A'\simeq\A'^0\oplus\A'^+$ is still
valid. 
\begin{lem}
Let $K=\epsh{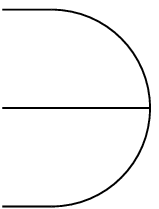}{2.5ex}\circ K'$ for some $K'\in\calD([0],[3])$ then
$$\Phi_\DD(K')=\Phi'_\DD(K).\Phi_\DD(\epsh{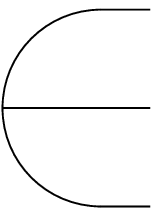}{3ex})$$
Furthermore, if $K\in\A([n])=\calD([0],[n])$ is sent on zero by
$\Phi_\DD$ then for any $K'\in\A'^+([n])\subset\calD'([n],[0])$, one has
$\Phi'_\DD(K'\circ K)=0$.
\end{lem}
Remark that the second assertion is false for $K'\in\A'^0([n])$!
\begin{theo} \label{skein}
Let $T=\epsh{triling}{2.5ex}\circ T'$ for some
$T'\in\Tp\big(\emptyset,(\circ\circ)\circ\big)$ then 
$$Z_{\DD,L}(T')=\wt Z(T).\Phi_\DD(\epsh{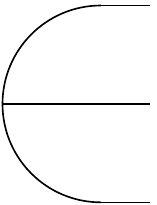}{3ex})$$ 
% %
Thus on proper knotted trivalent graph, $\wt Z$ can be computed using
$Z_{\DD,L}$. For links, one can compute the variation of $\wt Z$ when
one changes a crossing with:
\begin{equation*}
\wt Z\left(\epsh{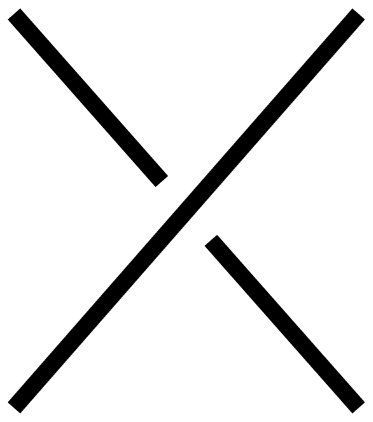}{4ex}-\epsh{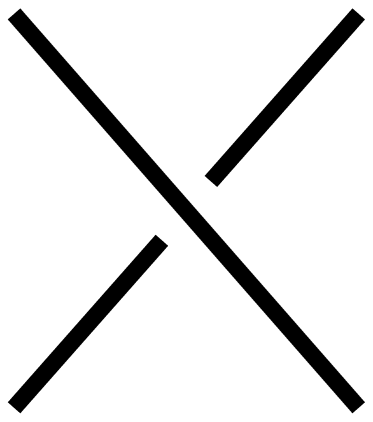}{4ex}\right)=-\frac12\wt
Z\left(\epsh{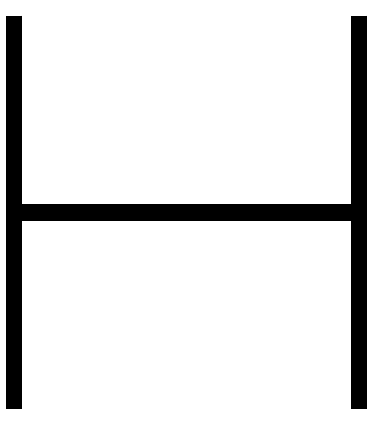}{4ex}\,-\epsh{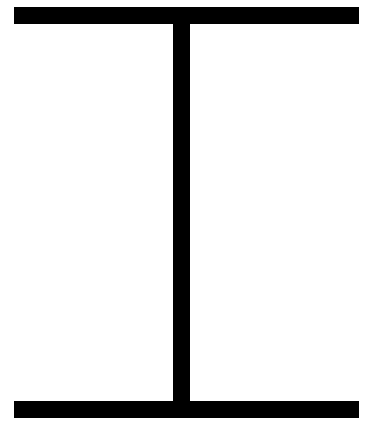}{4ex}+\frac12\left(\epsh{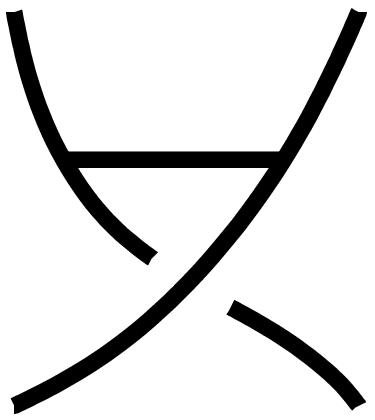}{4ex}+
\epsh{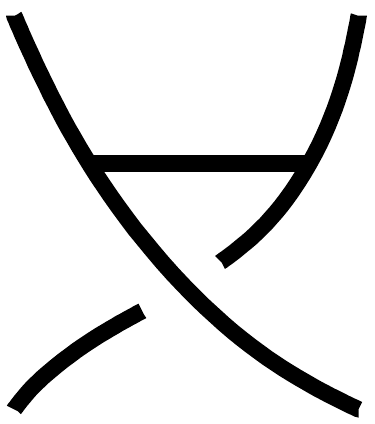}{4ex}\right)\right)
\end{equation*}
Furthermore, if $L_1$ and $L_2$ are links, $L_1'$ and $L_2'$ are
proper knotted trivalent graphs, 

$\wt Z(L_1\sqcup L_2)=\wt Z(L_1)+\wt Z(L_1)$, $\wt
Z(L_1'\sqcup L_2')=0$ and $\wt Z(L_1\sqcup L_1')=\wt Z(L_1')$.
\end{theo}

\begin{conj}\label{conjZ}
The value of $\wt Z$ on the unframed unknot is obtained by removing
the term with degree $-1$ in
$\frac{\sigma_++\sigma_--6}{4(\sigma_+-\sigma_-)}$ where $\sigma_\pm$
is defined as in \fullref{conjI} with $q_i=\e^{\frac{a_i}2}$.
\end{conj}

\section{The quantum group invariant}
%%%%%%%%%%%%%%%%%%%%%%%%%%%%%%%%%%%%%%%%%%%%%%%%%%%%%%%%%%%%%%%%%%%%%%%%%%%%%%
\subsection{The quantum group $\Uh$ and the functor $Q_{\DDa,L}$}
%%%%%%%%%%%%%%%%%%%%%%%%%%%%%%%%%%%%%%%%%%%%%%%%%%%%%%%%%%%%%%%%%%%%%%%%%%%%%%
\subsubsection{The quantum group $\Uh$}
Remark that there is three simple root systems for $\DDa$. Here and in
\fullref{DDa}, the presentations of the algebra are based on the
distinguish simple root system (with the smallest number of odd simple
roots) of $\DDa$. Unfortunately, this simple root system (contrary to
the simple root system with three odd simple roots) breaks the
symmetry that appear in \fullref{DD}. This symmetry, hidden in the
presentation of $\DDa$ \fullref{DDa} seem totally lost with the
deformation $\Uh$ of Y.M. Zou \cite{Zou}.

An universal $R$--matrix for $\DDa$ has been computed by H. Thys 
\cite{Th}. It depends of the three parameters $q_1$, $q_2$ and $q_3$
where $q_3=q_2^\alpha$ and $q_1q_2q_3=1$.

In the following, we will adopt the following notation for $i=1\cdots3$:
$$[n]_i=q_i^n-q_i^{-n}$$ 
% %
There exists an unique $17$--dimensional irreducible representation
$\rho$ of $\Uh$. Its classical limit is the adjoint representation of
the Lie superalgebras $\DDa$. This $\Uh$--module $L$ is autodual (there
is a (unique up to a scalar) bilinear map $\beta:L\otimes
L\go\CC[[h]]$), has the set of roots for set of weights, and has a
(unique up to a scalar) bilinear map $\gamma:L\otimes L\go L$ (whose
classical limit is the Lie bracket). 

We have computed using Maple
the $17^2\times17^2$ $R$--matrix for $L$, the tensor realizing the
duality $\beta$, its dual and the tensor $\gamma$. For a good choice
of a basis of $L$, and a good normalization of $\beta$ and $\gamma$,
all the coordinates of these tensors are in the ring
$\ZZ[\frac12,q_1,q_2,q_3,\Q41^{-1},\Q42^{-1},\Q43^{-1}]/(q_1q_2q_3=1)$.
The computations with Maple are available on the author's web page.
% %

%%%%%%%%%%%%%%%%%%%%%%%%%%%%%%%%%%%%%%%%%%%%%%%%%%%%%%%%%%%%%%%%%%%%%%%%%%%%%%
\subsubsection{$Q_{\DDa,L}$ and the quantum Jacobi relation}\label{QD}
As usual one has a functor $Q_{\DDa,L}$ from $\T$ to $\Mod_{\Uh}$
sending $[1]\in\Obj(\T)$ to the $\Uh$--module $L$. This givea an
invariant of framed knotted trivalent graphs with values in
$$\ZZ[\frac12,q_1,q_2,q_3,\Q41^{-1},\Q42^{-1},\Q43^{-1}]/(q_1q_2q_3=1).$$

\begin{theo} \label{relQ}
The following elements are in the kernel of $Q_{\DDa,L}$:
$$\epsh{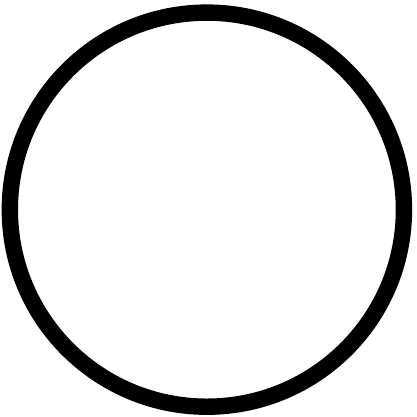}{4ex}-1\ets
\epsh{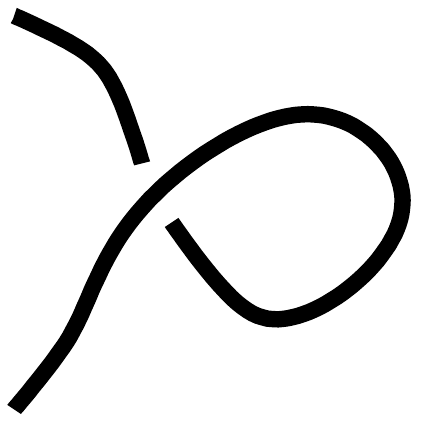}{4ex}-\epsh{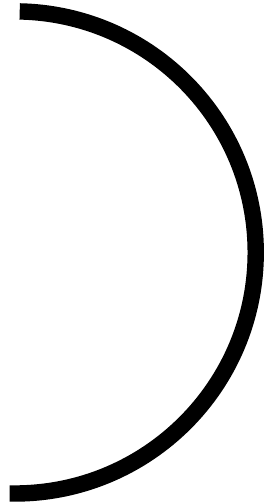}{4ex}\ets
\epsh{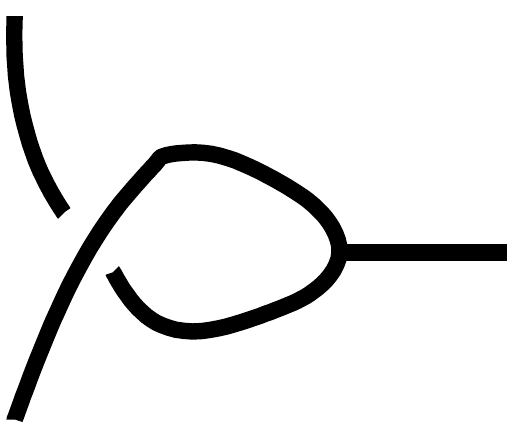}{4ex}+\epsh{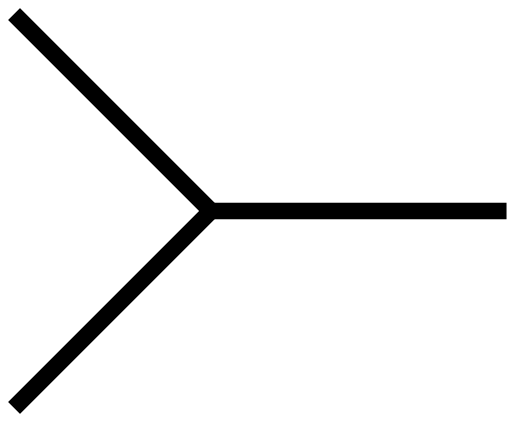}{4ex}\et \epsh{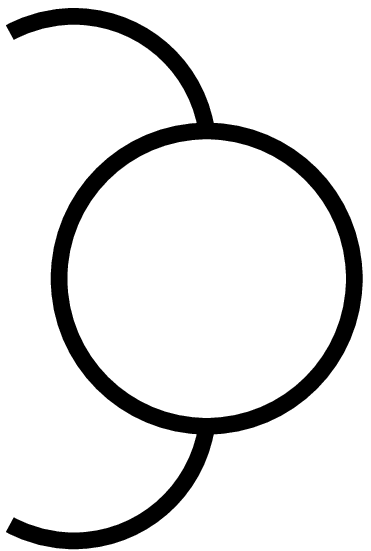}{4ex}$$
% twist-1 , Cr@Twist+Cr , Bulle 
Furthermore, $\Mod_{\Uh}(L^{\otimes 2},L^{\otimes 2})$ has dimension $6$
and is generated by the images by $Q_{\DDa,L}$ of the powers of the
half twist. The ``quantum Jacobi relation'' is satisfied by
$Q_{\DDa,L}$:
$$\epsh{tangTp}{4ex}-\epsh{tangTm}{4ex}\stackrel{Q_{\DDa,L}}=\phi\left(
\epsh{tangH}{4ex}-\epsh{tangI}{4ex}+\frac12\left(\epsh{tangXp}{4ex}+
\epsh{tangXm}{4ex}\right)\right)$$ 
% %
where $\phi$ depends of the normalization chosen for $\beta:L\otimes
L\go\CC[[h]]$ and $\gamma:L\otimes L\go L$.
\end{theo}
\begin{proof}
This is a computation made with Maple.
\end{proof}

\begin{proof}[Proof of \fullref{Q=0}]
Every proper knotted trivalent graph $T$ can be written $T_1\circ T_2$ with
$T_1\in\T(0,3)$ and $T_2 \in\T(3,0)$. But the spaces
$\Mod_\Uh(L^{\otimes p},L^{\otimes q})$ with $p+q=3$ are all
isomorphic with dimension $1$ so $Q_{\DDa,L}(T)$ is proportional to
the image by $Q_{\DDa,L}$ of the knotted trivalent graph $\Theta$ which is $0$. As a
consequence, the ``quantum Jacobi relation'' implies that $Q_{\DDa,L}$
is unchanged when one changes the crossings of a link and so is constant
equal to its value on the unlink which is $1$.
\end{proof}

We will need the following:
\begin{coro} \label{mutat} If $r_\pi:\T(2,2)\go\T(2,2)$ is the map
  induced by the rotation around the line
  $\{\frac12\}\times\{\frac32\}\times\RR$ by $\pi$ then one has
  $Q_{\DDa,L}\circ r_\pi=Q_{\DDa,L}$ on $\T(2,2)$.
\end{coro}

%%%%%%%%%%%%%%%%%%%%%%%%%%%%%%%%%%%%%%%%%%%%%%%%%%%%%%%%%%%%%
\subsection{Renormalization of $Q_{\DDa,L}$}
The idea is to define a renormalization $\wt Q$ of $Q_{\DDa,L}$ using
some relation between the two as in \fullref{skein}. The
demonstration of the invariance is then not trivial{\footnote{We remark
that according to the new results of \cite{Geer2}, the invariance of
$\wt Z$ implies the existence and the invariance of $\wt Q$}} but it
is made by analogy to some demonstrations for weight functions. We
give the steps of the demonstration:
\begin{itemize}
\item $\wt Q$ is well defined on proper knotted trivalent graphs.
\item Using the quantum Jacobi relation, it can be extended to an
  invariant of singular links with one double point.
\item This invariant can be integrated to a link invariant $\wbtot Q$.
\end{itemize}

\subsubsection{$\wt Q$ for Proper knotted trivalent graphs and Singular Link}
\begin{theo}
Let $L=\epsh{triling}{2.5ex}\circ T$ for some $T\in\T(0,3)$ then the
scalar $\wt Q(L)$ defined by $Q_{\DDa,L}(T)=\wt
Q(L).Q_{\DDa,L}(\epsh{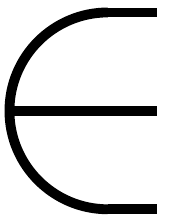}{3ex})$ is independent of $T$.
\end{theo}
\begin{proof}
First remark that by \fullref{relQ}, $Q_{\DDa,L}(T)$ does not
depend of the framing of $T$ and that the braid group
$B_3\subset\T(3,3)$ acts on $Q_{\DDa,L}(\T(0,3))$ by multiplication by
the signature (a braid $b$ with projection $\sigma\in\S_3$ act by the
multiplication by the signature of $\sigma\in\S_3$ (cf the third
relation of \fullref{relQ})). So it is easily seen that
$Q_{\DDa,L}(T)$ depends only of the choice of the trivalent vertex of
$L$ that is removed in $T$.

Second, if $L$ is a disjoint union of knotted trivalent graph then clearly,
$Q_{\DDa,L}(T)=0$.

Third, by \fullref{mutat} applied to a tangle $P\in\T(2,2)$, one
has:
$$Q_{\DDa,L}\left(\epsh{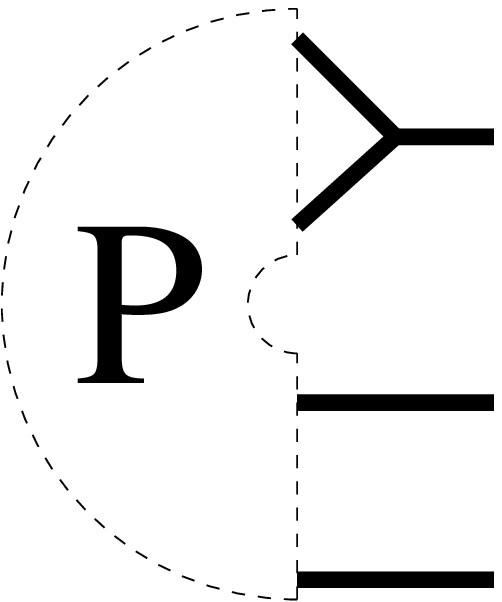}{4ex}\right)=
Q_{\DDa,L}\left(\epsh{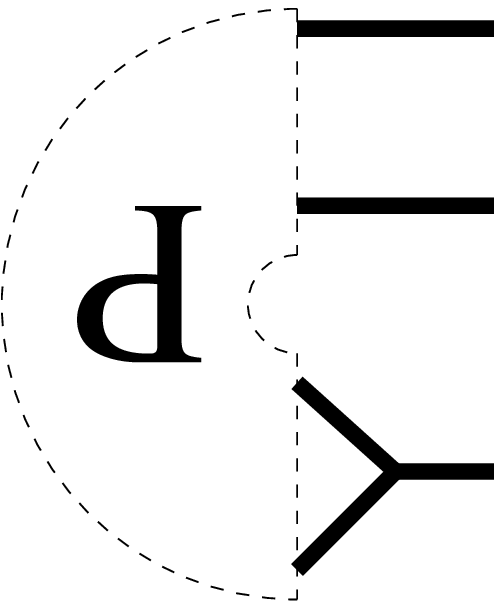}{4ex}\right)=
Q_{\DDa,L}\left(\epsh{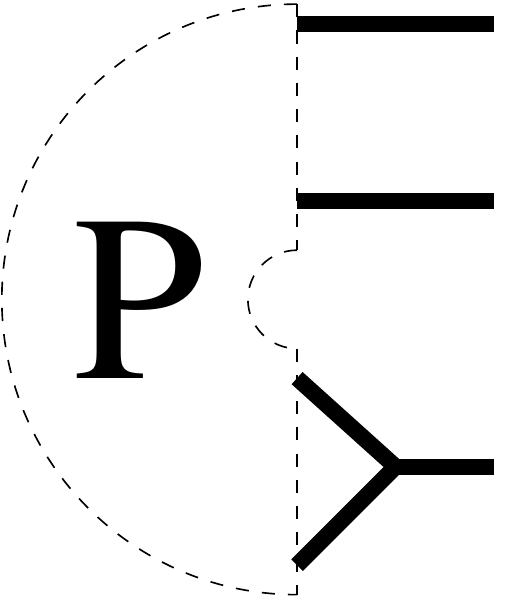}{4ex}\right)$$ 
% %
So $\wt Q(L)$ is unchanged when one chooses any $3$--vertex in the same
connected component of $L$ and the theorem is proved for connected
knotted trivalent graph.

Last, consider two trivalent tangles $T$ and $T'$ in $\T(0,3)$ giving
the same knotted trivalent graph. One can find a trivalent tangle $\wb
T\in\T(0,6)$ such that
$$T=(\Id^{\otimes 3}\otimes\epsh{triling}{2.5ex})\circ\wb T\textrm{
  and }T'=(\epsh{triling}{2.5ex}\otimes\Id^{\otimes 3})\circ\wb T$$
% %
Now one can use the quantum Jacobi relation to change the crossing in
$\wb T$ and find a sum of trivalent tangles $\wt T$ such that
$Q_{\DDa,L}(\wb T)=Q_{\DDa,L}(\wt T)$ and the trivalent tangles that appear
in $\wt T$ are either the tensor product of two trivalent tangles in
$\T(0,3)$ (so do not contribute to $Q_{\DDa,L}(T)$ nor to
$Q_{\DDa,L}(T')$) or are trivalent tangles with at least one component
intersecting both the sets of univalent vertices $\{1,2,3\}$ and
$\{4,5,6\}$ (so they contribute for the same as they give connected
knotted trivalent graph).
\end{proof}
\begin{Rq}
By \fullref{relQ}, $\wt Q$ is independent of the framing.
\end{Rq}
\begin{defi}\label{Qsing}
If $L$ is a framed oriented link with one double point, we define $\wt
Q(L)$ by the following substitution:
 \begin{equation*}
\wbtot Q\left(\epsh{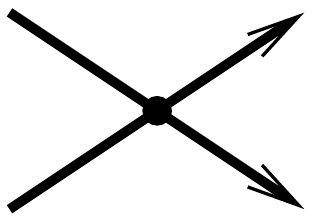}{4ex}\right)=\wt
Q\left(\epsh{tangH}{4ex}-\epsh{tangI}{4ex}+\frac12\left(\epsh{tangXp}{4ex}+
\epsh{tangXm}{4ex}\right)\right)
\end{equation*}
\end{defi}
Remark that the orientation of $L$ is forgotten in the right hand side.

\begin{prop} For framed oriented links with one double point as
  follow, one has 
$$\wbtot Q\left(\epsh{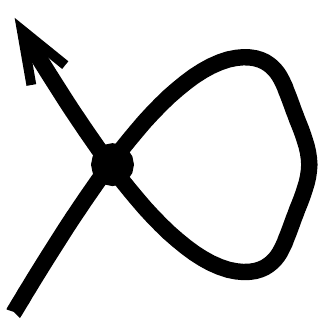}{4ex}\right)=2$$
\end{prop}

For a framed link $L$, let $w(L)$ denotes the diagonal writhe of $L$
(ie, the trace of the linking matrix of any orientation of $L$). Then
$w$ extends to an invariant of framed oriented links with one double
point which also satisfies:
$$w\left(\epsh{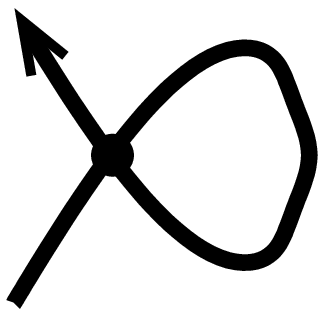}{4ex}\right)=2$$

\subsubsection{Integration of $\wbtot Q$}\label{GQD}
In \cite[Theorem 1]{St} T Stanford gives local conditions for a
singular link invariant to be integrable to a link invariant:
specialized in our context, it gives:
\begin{theo}[Stanford]\label{Sth} Let $\LL^{(1)}$ be the set of
  isotopy classes of singular links in $\RR^3$ with one double point
  and let $k$ be a field of characteristic zero.

Then, for any
  finite type singular link invariant $f:\LL^{(1)}\go k$, there exists
  a link invariant $F:\LL\go R$, such that
$$f\left(\epsh{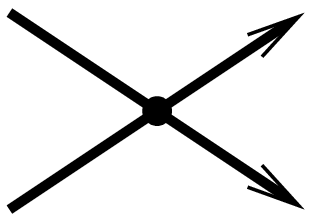}{4ex}\right)=
  F\left(\epsh{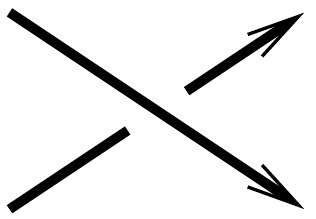}{4ex}\right)-F\left(\epsh{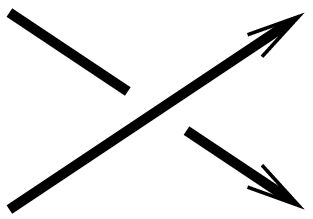}{4ex}\right)$$
iff
\begin{equation}\label{cocycle}
f\left(\epsh{tangFrams}{4ex}\right)=0\et
f(L_{\times\,+})-f(L_{\times\,-})=f(L_{+\,\times})-f(L_{-\,\times})
\end{equation}
(where $L_{*\,*}$ denotes some desingularisations of a any singular link
$L_{\times\times}$ with two double points).
\end{theo}

\begin{theo}
$\wbtot Q-w$ satisfies the two conditions \ref{cocycle} of 
  \fullref{Sth} and so $\wbtot Q$ extends in an unique way to a framed link
  invariant which takes value zero on the unlink.

Furthermore, $\wbtot Q$ takes value in the ring
$\QQ[q_1,q_2,q_3,\Q41^{-1},\Q42^{-1},\Q43^{-1}]_{/q_1q_2q_3=1}$ and
satisfy $\wbtot Q(L\#L')=\wbtot Q(L)+\wbtot Q(L')$ where $L\#L'$
denotes a connected sum of $L$ and $L'$ along one of their components.
\end{theo}
\begin{proof}
For a singular link $L_{\times\times}$ with two double points, the two
terms $\wbtot Q(L_{\times\,+}-L_{\times\,-})$ and $\wbtot
Q(L_{+\,\times}-L_{-\,\times})$ are equal to $\wbtot Q(K')$ where $K'$
is the sum of knotted trivalent graphs obtained by replacing the two
singular points of $L_{\times\times}$ as in \fullref{Qsing}.
\end{proof}

\begin{conj}Relation between $\wt Z$ and $\wt Q$\footnote{The first
    part of this conjecture is now proven by N Geer in \cite{Geer2}.}
    \label{conjQ}
%\begin{itemlist}
\begin{enumerate}
\item For any proper knotted trivalent graph $L$ one has $\wt Q(L)=\wt
Z(L)$.
\item For any framed link $L$ with $n$ components,\footnote{after
  removing the term of degree $-1$} 
$$\wbtot Q(L)=2\wt
  Z(L)-n\frac{\sigma_++\sigma_--6}{2(\sigma_+-\sigma_-)}$$
\item $\wbtot Q$  takes value in the polynomial algebra
$\ZZ[\sigma_+,\sigma_-]$
where $\sigma_+=(q_1^2+q_2^2+q_3^2)$ and
  $\sigma_-=(q_1^{-2}+q_2^{-2}+q_3^{-2})$.
%\end{itemlist}
\end{enumerate}
\end{conj}
Remark that $1$ implies $2$ (with \fullref{conjZ}), and
the fact that the values of $\wbtot Q$ are symetrics in the three
variables.

%%%%%%%%%%%%%%%%%%%%%%%%%%%%%%%%%%%%%%%%%%%%%%%%%%%%%%%%%%%%%
\section{Properties of the invariants}
\subsection{The common specialisation with the Kauffman polynomial}
Remember that the specialisation $\alpha\in\{-2,-\frac12,1\}$ of $\DDa$
give a Lie superalgebra isomorphic with $\osp(4,2)$. So in this case,
$\Uh$ admit a $6$--dimensional representation which satisfies the
skein relations of the Kauffman polynomial:
$$K\left(\epsh{tangTp}{4ex}-\epsh{tangTm}{4ex}\right)=
(s-s^{-1})K\left(\epsh{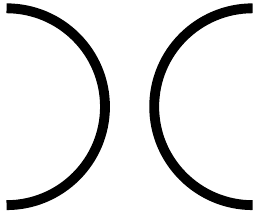}{4ex}-\epsh{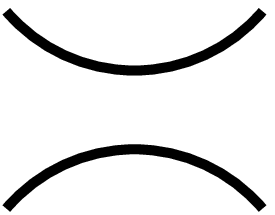}{4ex}\right)\qua
\text{and}\qua K\left(\epsh{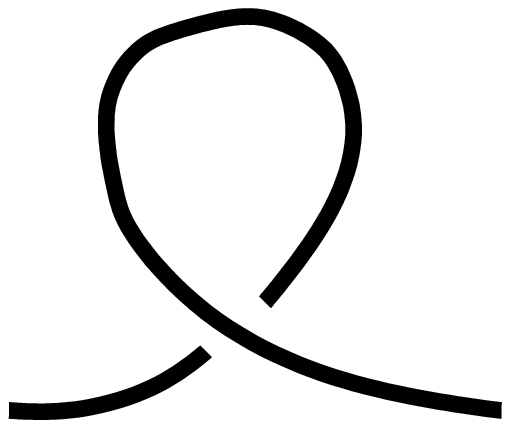}{4ex}\right)=\alpha
K\left(\epsh{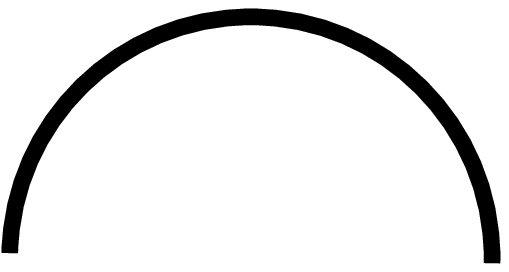}{3ex}\right)$$ for $\alpha=s$.

As these skein relations determine the tangle invariant, the
specialisations of the functor $Z_\DD$ obtained by setting
$a_1=a_2=a$, $a_3=-2a$ (or any permutation of this) and the
specialisations of the functor $Q_{\DDa,L}$ obtained by setting
$\alpha\in\{-2,-\frac12,1\}$ (ie, $q_1=q_2=s^{-1}$, $q_3=s^{2}$) are
both equivalent to the $\alpha=s$ specialisation of the ``adjoint''
Kauffman skein quotient which is obtained by cabling each component of
a tangle with the following projector of $\T([2],[2])$:
$$\frac1{s+s^{-1}}\left(s\epsh{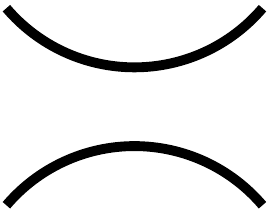}{4ex}-\epsh{tangTm}{4ex}-
\frac{s-s^{-1}}{\alpha s^{-1}+1}\epsh{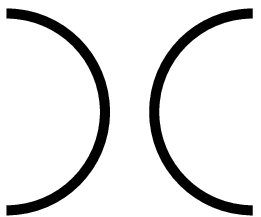}{4ex}\right)$$ and
imposing the Kauffman skein relations.

Let $K_\ad$ be the framed link invariant obtained by cabling each
component of a framed link with the previous projector and computing
its Kauffman polynomial then:
\begin{theo} \label{Kauf} Let $\phi$ be the specialisation
  $\phi(q_1)=\phi(q_2)=s^{-1}$, $\phi(q_3)=s^2$ (so that
  $\phi(\sigma_+)=2s^{-2}+s^{4}$ and
  $\phi(\sigma_-)=2s^2+s^{-4}$). Then for any framed link $L$,
  $\phi(\wbtot Q(L))$ and $\phi(\wt Z(L))$ are related as in
  \fullref{conjQ} and 
\begin{gather*}
K_\ad(L)|_{\alpha=s}=1\\
\left.\left(\frac
 {K_\ad(L)-1}{\alpha-s}\right)\right|_{\alpha=s}=\frac{2}{s}\phi(\wbtot
 Q(L))
\end{gather*}
\end{theo}

%%%%%%%%%%%%%%%%%%%%%%%%%%%%%%%%%%%%%%%%%%%%%%%%%%%%%%%%%%%%%
\subsection{The common specialisation with the HOMFLY-PT
 polynomial}\label{Homflypt} 
% %
It would be more difficult to make appear the common specialisation of
$\wbtot Q$ with the HOMFLY-PT polynomial. This should appear for the
degenerate specialisation $\alpha\in\{0,-1,+\infty\}$ of $\DDa$. We
state the existing relation between $\wt Z$ and HOMFLY-PT and we just
state a conjecture for the relation between $\wbtot Q$ and
HOMFLY-PT\footnote{This conjecture now follows from the work of N Geer
(see \cite{Geer2}).}.

The HOMFLY-PT polynomial of an oriented link $L\neq \emptyset$ is an
element $P(v,z)\in\ZZ[v^{\pm},z^{\pm}]$ which is equal to $1$
  for the unknot and satisfy the skein relation:
$$v^{-1}P\left(\epsh{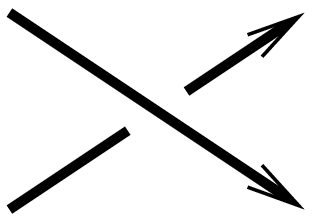}{4ex}\right)-
  vP\left(\epsh{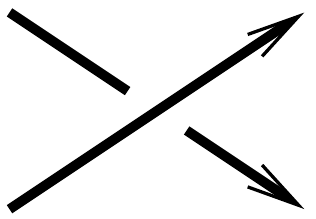}{4ex}\right)=
  zP\left(\epsh{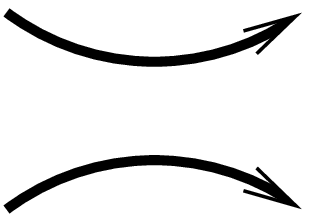}{4ex}\right)$$ 
% %
If $W$ denotes the total writhe of an oriented framed link $L$
(ie, the total algebraic number of crossing of $L$) then we get an
oriented framed link invariant of $L$:
$H(\lambda,v,z)\in\ZZ[\lambda^{\pm},v^{\pm},z^{\pm}]$ by the formula
$$H(L)=\left|\begin{array}{l}1\text{ if }L=\emptyset
\lambda^{W(L)}P(L)\times\frac{v^{-1}-v}z\text{
  else}\end{array}\right.$$
$H$ satisfy the skein relations:
$$(v\lambda)^{-1}H\left(\epsh{tanghTpo}{4ex}\right)- (v\lambda)
  H\left(\epsh{tanghTmo}{4ex}\right)=
  zH\left(\epsh{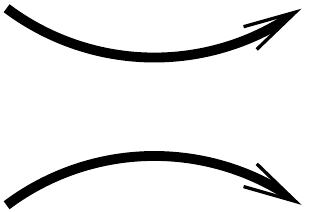}{4ex}\right)\text{ and }
  H\left(\epsh{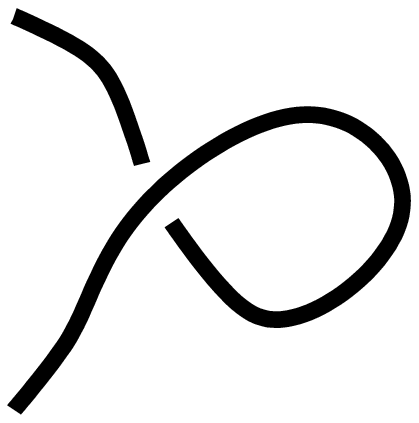}{4ex}\right)=\lambda
  H\left(\epsh{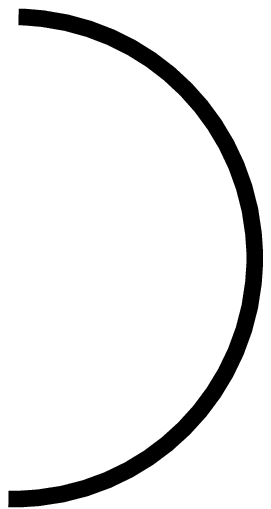}{4ex}\right)$$ 
% %
We define the adjoint HOMFLY-PT polynomial $H_\ad$ of a framed
unoriented link $L$ as the $H$ polynomial of the framed oriented link
obtained by cabling each component of $L$ with the following:
$$\left(\,{\begin{array}{c} \hspace{-1.3mm}
        \raisebox{2pt}{\includegraphics[width=5ex]{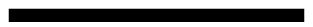}}
        \hspace{-1.9mm}\end{array}}\,\right)\mapsto
  \left(\epsw{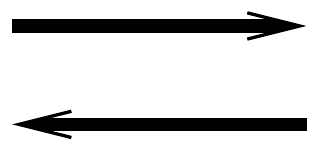}{5ex}- \emptyset\right)$$ 
% %
Remark that the cabled link has total writhe $0$, so $H_\ad(L)$ is
proportional to $P(v,z)$ of the cabled link on $L$ and so lies in
$\ZZ[v^{\pm},z^{\pm}]$.  One can also compute $H_\ad(q^n,q^{-1}-q)$ by
the way of a quantum group $U_q(\sll_n)$ and its ``adjoint''
$(n^2-1)$--dimensional representation (here ``adjoint'' mean the
quantum analogue of the adjoint representation of $\sll_n$) or
equivalently by composing $Z_\ad$ with the $\sll_n$ weight system.

In fact $f\circ H=\Phi_\sll\circ Z$ where the weight function
$\Phi_\sll$ takes values in the ring $\QQ[\delta^{\pm},h]$ and $f$ is
the ring morphism such that $f(v)=\e^{-\frac h2{\delta}}$,
$f(z)=2\sinh(h/2)=\e^{\frac12 h}-\e^{-\frac12h}$ and
$f(\lambda)=\e^{\frac h2(\delta-\frac1\delta)}$. Remark that there
exists a character $\chi_\sll$ on $\Lambda$ associated with
$\Phi_\sll$ and whose values belong to $\QQ[[\delta h,h^2]]$.

We show in \cite{moi} that on $\Lambda$, the map $\chi_\sll$ modulo
$\delta$ and $\chi_\DD$ modulo $\sigma_3$ were both zero and that
$\chi_\sll$ modulo $\delta^2$ and $\chi_\DD$ modulo $\sigma_3^2$ were
equal up to renormalization:

For $\lambda\in\Lambda$ of degree $2p+1\geq3$, if $\chi_{\DD}(\lambda)=
\mu(\lambda) \sigma_3(\sigma_2)^{p-1} +O(\sigma_3^2)$

then $\Phi_\sll(\lambda)=(-1)^p\mu(\lambda)\delta h^{2p}
+O(\delta^2)$. (Here and after, $O(x)$ denotes an element of the ideal
generated by $x$).

We call $\psi$ the specialization defined by $\psi(\sigma_3)=-\delta
h^2$, $\psi(\sigma_2)=-h^2$ (So that $\psi(\sigma_\pm)=
1+\e^h+\e^{-h}\mp\frac\delta2(\e^h+\e^{-h}-2)+O(\delta^2)=
f(z^2+3\mp\frac\delta2z^2)+O(\delta^2)$).

For a knot $K$ closure of $T\in\T([1],[1])$, we have

$Z_\ad(T)=1+\frac{w(K)}2\epsh{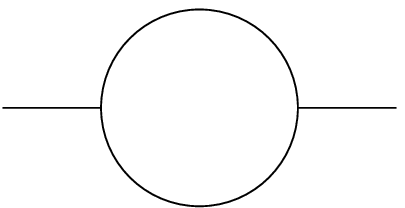}{2ex}
+t(\frac{w(K)}2)^2\epsh{geneF2}{2ex}-V_2(K)t\epsh{geneF2}{2ex}
+\lambda.\epsh{geneF2}{2ex}\in\calD([1],[1])$

where $w(K)$ is the writhe of $K$, $V_2(K)$ is the (standardly
normalized) type $2$ Vassiliev invariant of the knot $K$ (ie, the
coefficient of $z^2$ in $P(1,z)$) and $\lambda\in\Lambda$, $\lambda$
nul in degree $\leq1$.

So $\Phi_\sll(Z_\ad(T))= 1+ w(K)\delta
h+\left(\frac{w(K)^2}2-2V_2(K)\right)\delta^2h^2+
2\delta\chi_\sll(\lambda)$

whereas $\wt Z(K)=\wt Z(U_0)+\frac{w(K)}2+\chi_{\DD}(\lambda)$ (with
$U_0$ the unknot).

So we get
% using $v-v^{-1}=-\delta+O(\delta^3)$ 
$$\frac{f(H_\ad(K))}{f(H_\ad(U_0))}=1
+\delta^2\left(\frac{w(K)^2}2-2V_2(K)\right) +2\delta\psi(\wt
Z(K)-\wt Z(U_0))+O(\delta^3)$$

So using 1) and 3) of \fullref{conjQ} one has
\begin{conj} For $K$ a $0$--framed knot,\footnote{This conjecture
now follows from the work of N Geer (see \cite{Geer2}).}
 $$ \left.\frac{\frac{f(H_\ad(K))}{f(H_\ad(U_0))}-1}
 {\left(v-\frac1v\right)^2}\right|_{v=1}= -2V_2(K)-\frac1{z^2}
 \left.\frac{\wbtot Q(K)}{\sigma_+-\sigma_-}
 \right|_{\sigma_+=\sigma_-=z^2+3}$$
\end{conj}

%%%%%%%%%%%%%%%%%%%%%%%%%%%%%%%%%%%%%%%%%%%%%%%%%%%%%%%%%%%%%
\subsection{Example of computation}
With some computations on Maple, we found that in the base of
$\Mod_\Uh(L^{\otimes 2},L^{\otimes 2})$ given by
$(U,T^{-2},T^{-1},Id,T,T^2)$ where
$T=Q_{\DDa,L}\left(\epsh{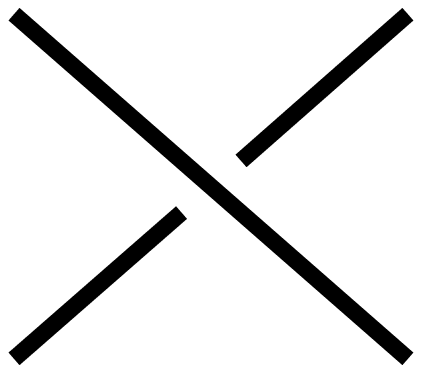}{3ex}\right)$ is the positive half
twist and $U=Q_{\DDa,L}\left(\epsh{tanghU}{3ex}\right)$, one has 
\begin{equation}\label{rotT2}
Q_{\DDa,L}\left(\epsh{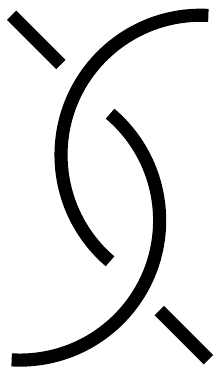}{5ex}\right)=\left [\begin
    {array}{c}1-2(\sigma_+-\sigma_-)\\ 0\\ -2+\sigma_-\\
    -1+2(\sigma_+-\sigma_-)\\ 2-\sigma_-\\ 1\end{array}\right ]
\end{equation}
% %
Furthermore,
\begin{equation}\label{T3}
Q_{\DDa,L}\left(T^3\right)=
\left[\begin{array}{c}4 \left (\sigma_+-\sigma_-\right )\\ 1\\
2-\sigma_+\\ 1-2\sigma_++\sigma_-\\-1-\sigma_++2\sigma_-\\-2+\sigma_-
\end{array}\right]
\end{equation}
Let $K_{2n+1}$ be the knot obtained as the closure of $T^{2n+1}$ and
$I_n$ be the knotted trivalent graph obtained as the closure of
$I\circ T^n$ where $I=\epsh{tangI}{2ex}$.

So by \fullref{QbKnot} one has:
\begin{align*}
&\wbtot Q(K_{2n+1})-\wbtot Q(K_{2n-1})=\left((-1)^{2n}-\wt
Q(I_{2n})+\frac12(-\wt Q(I_{2n-1})-\wt Q(I_{2n+1}))\right)\\
&\qquad=1-\frac12\left(\wt Q(I_{2n-1})+2\wt
Q(I_{2n})+\wt Q(I_{2n+1})\right)\\
&\text{with }\\
&\wbtot Q(K_1)=1=-\wbtot Q(K_{-1}),\quad\wt Q(I_0)=0,\quad \wt
Q(I_{\pm1})=-1\\ 
&\wt Q(I_{2})=2(\sigma_+-\sigma_-),\quad \wt
Q(I_{-2})=-2(\sigma_+-\sigma_-)\quad \text{by equation
  \ref{rotT2}}\\ 
&\wt Q(I_{n+3})=4(\sigma_+-\sigma_-)+\wt Q(I_{n-2})+(2-\sigma_+)\wt
Q(I_{n-1})+(1-2\sigma_++\sigma_-)\wt
Q(I_{n})\\
&\qquad+(-1-\sigma_++2\sigma_-)\wt Q(I_{n+1})+(-2+\sigma_-)\wt
Q(I_{n+2})\quad\text{by equation \ref{T3}}
\end{align*}
Thus one can compute:
\begin{align*}
\wbtot Q(K_{3})=&\,3-(\sigma_+-\sigma_-)(2 +\sigma_-)\\
\wbtot Q(K_{5})=&\,5+(\sigma_+-\sigma_-)(-6 +2\sigma_+ -4\sigma_-
+2\sigma_+\sigma_- -2\sigma_-^2 -\sigma_-^3)\\
\vdots&
\end{align*}
The same method gives $\wbtot Q($Hopf link$)=\sigma_--\sigma_+$.\\
And to compare with \fullref{Kauf},
(with $U_0$ the unframed unknot) we see that:
\begin{align*}
K_\ad(U_0)&={\frac {\left (\alpha^2-1\right ) \left ({s}^{3}+
    \alpha\right ) \left (s\alpha-1\right
)s}{{\alpha}^{2}\left ({s}^{4}-1 \right )\left (s^2-1\right )}}\\
% %
=1+&{\frac {\left (\alpha-s\right )\left
({\alpha}^{3}{s}^{2}+{\alpha}^{2}({s}^{5}+{s}^{3}-s)
+\alpha({s}^{4}-{s}^{2}-1) -{s}^{3}\right )}{{\alpha}^{2}\left
({s}^{4}-1 \right )\left (s^2-1\right )}}\\
\left.\frac{K_\ad(U_0)-1}{\alpha-s}\right|_{\alpha=s}&={\frac
{{s}^{4}+4\,{s}^{2}+1}{s\left (s^4-1\right )}}\\
%%%
% %
\frac{K_\ad(K_3)}{K_\ad(U_0)}=(\alpha^2-s^2)&\left({\frac
{{s}^{12}+{s}^{8}+{s}^{6}+1}{{s}^{10}}}+
% %
{\frac {\left (s^4-1\right )\left ({s}^{6}+1\right
)}{{s}^{7}\alpha}}\right.\\
% %
&\left.-{\frac {{s}^{12}-{s}^{10}-{s}^{8}+2\,{s}^{6}
    -{s}^{2}+1}{{s}^{6}{\alpha}^{2}}}-{\frac {\left (s^4-1\right
    )\left ({s}^{6}-{s}^{2}+1 \right )}{{s}^{3}{\alpha}^{3}}}\right.\\
% %
&\left.-{\frac {\left ({s}^{4}-1\right )\left (s^2-1\right
)}{{\alpha}^{4}}}
% %
\right)
% %
\end{align*}

\begin{align*}
\left.\frac{K_\ad(K_3)-1}{\alpha-s}\right|_{\alpha=s}
&=\frac2s\left(\phi(\wbtot Q(K_3))+ {\frac {{s}^{4}+4{s}^{2}+1}{s\left
(s^4-1\right )}}\right)
\end{align*}
And to compare with \fullref{Homflypt}
\begin{align*}
H_\ad(U_0)&={\frac { \left( {v}^{2}+zv-1 \right) \left( {v}^{2}-zv-1
\right) }{{z }^{2}{v}^{2}}}=-1+\frac{\left( v-\frac1v \right)
^{2}}{z^2} \\ 
% %
\frac{H_\ad(K_3)}{H_\ad(U_0)}&=1-3\left( v-\frac1v \right)+\left(
v-\frac1v \right) ^{2}\left({\frac {v+4}{v+1}+ \left( {v}^{2}+4
\right) {z}^{2}+{z}^{4}} \right)\\
% %
f\left(\frac{H_\ad(K_3)}{H_\ad(U_0)}\right)&=1+3\delta+ \left(\frac52
+5{z}^{2}+ {z}^{4}\right) {\delta}^{2} +O(\delta^3)\\
% %
&=1+3\delta+ \delta^2\left(\frac92
-2\right)+\delta^2z^2\Big(2+(z^2+3)\Big) +O(\delta^3)\\
% %
&=1+\delta^2\left(\frac{W(K_3)}2-V_2(K_3)\right)+\delta\psi(\wbtot
Q(K_3))+O(\delta^3) 
\end{align*}

\bibliographystyle{gtart} \bibliography{link}

\end{document}